\documentclass[12pt]{article}
\usepackage{amsthm,amsmath,amssymb,mathtools}
\usepackage{graphicx,psfrag,epsf}
\usepackage{enumerate}
\usepackage{url} % not crucial - just used below for the URL
\usepackage{latexsym}
\usepackage{color}

\definecolor{blue}{rgb}{0,0,0.9}
\definecolor{red}{rgb}{0.9,0,0}
\definecolor{tgreen}{rgb}{0,0.50,0.10}
\definecolor{violet}{rgb}{0.5804,0.0000,0.8275}

\newcommand{\blind}{0}

\newtheorem{theorem}{Theorem}

\newcommand{\E}{{\bf E}}
\newcommand{\Prob}{{\bf P}}
\def\QED{{\setlength{\fboxsep}{0pt}\setlength{\fboxrule}
{0.2pt}\fbox{\rule[0pt]{0pt}{1.3ex}\rule[0pt]{1.3ex}{0pt}}}}

 \usepackage{setspace}

\begin{document}

%%%%%%%%%%%%%%%%%%%%%%%%%%%%%%%%%%%%%%%%%%%%%%%%%%%%%%%%%%%%%%%%%%%%%%%%%%%%%%

\if0\blind
{
  \title{\bf An elementary approach for minimax estimation of Bernoulli proportion in the restricted parameter space}
 \author{Heejune Sheen\thanks{ ({\tt brianshn@gatech.edu})}
%    The authors gratefully acknowledge \textit{please remember to list all relevant funding sources in the unblinded version}}\hspace{.2cm}\\
    %School of Industrial and Systems Engineering, \\ Georgia Institute of Technology, Atlanta, GA, USA   \\
  %           ({\tt brianshn@gatech.edu}). \\
 \   and
   Yajun Mei\thanks{({\tt ymei@isye.gatech.edu})} \\~\\
  H. Milton Stewart School of Industrial and Systems Engineering, \\ Georgia Institute of Technology, Atlanta, GA, USA    \\
 }
  \maketitle
} \fi

\if1\blind
{
  \bigskip
  \bigskip
  \bigskip
  \begin{center}
    {\LARGE\bf An elementary approach for minimax Estimation of Bernoulli proportion in the restricted parameter space}
    \author{}
%  \author{Heejune Sheen\thanks{
%    %The authors gratefully acknowledge \textit{please remember to list all relevant funding sources in the unblinded version}}\hspace{.2cm}\\
%    School of Industrial and Systems Engineering, Georgia Institute of Technology, Atlanta, GA, USA
%             ({\tt brianshn@gatech.edu}).\\
%    and \\
%   Yajun Mei\thanks{
%  School of Industrial and Systems Engineering, Georgia Institute of Technology, Atlanta, GA, USA
%             ({\tt yajun.mei@isye.gatech.edu}).}
%}
\end{center}
  \medskip
} \fi

\bigskip
\begin{abstract}
We present an elementary mathematical method to find the minimax estimator of the Bernoulli proportion $\theta$  under the squared error loss
when $\theta$ belongs to the restricted parameter space of the form $\Omega=[0, \eta]$ for some pre-specified constant $0\le \eta \le 1.$
This problem is inspired from the problem of estimating the rate of positive COVID-19 tests. The presented results and applications would be useful materials for both instructors and students when teaching point estimation in statistical or machine learning courses.
\end{abstract}


\vfill

\newpage
%\spacingset{1.25} % DON'T change the spacing!

\section{Introduction}
\label{sec:intro}

%Point estimation of model parameters have been an important topic in statistics, machine learning, and data science. Besides the maximum likelihood estimator (MLE), the method of moment (MOM), and Bayesian method, another important approach is the minimax estimator that minimizes the maximum risk.

Teaching the topic of minimax estimators is intriguing in the undergraduate or graduate level statistical inference course. On the one hand,
the applications of minimax estimators can be found in many fields such as
in statistics \cite{malinovsky2015note,yaacoub2018optimal,zinodiny2011bayes}, machine learning \cite{ben2007blind}, physics
\cite{ng2012simple,ng2012minimax} and finance \cite{chamberlain2000econometric}. On the other hand, the classical statistical textbook often only include the minimax estimators of Binomial or normal distributions that are derived through the Bayes estimators with the least favorable prior. There are also other approaches to prove minimax properties, such as information inequality \cite{Tsybakov:2008:INE:1522486} or the invariance methods \cite{kiefer1957invariance}, but they are generally too advanced to undergraduate or graduate students whose major might not be on statistics, e.g., those in machine learning, optimization, operation research, and data science. % However, in general it is challenging to derive the minimax estimators, and the standard statistical textbook provides very limited examples of minimax procedures.
This severely limits feasible practical exercises on the minimax estimators, and thus severely affects the interests of students with diverse backgrounds.

The goal of this note is to present an elementary mathematical method to find the minimax estimator of the Bernoulli proportion $\theta$  under the squared error loss when $\theta$ belongs to the restricted parameter space, thereby enriching teaching materials for both instructors and students when teaching minimax estimators in statistical or machine learning courses. Recall that the estimation problem including the minimax estimator on the proportion $\theta$ of the Bernoulli or binomial distributions has been well studied for the unrestricted parameter space in the literature. To be more specific,
%One important question on the minimax estimator is to estimate the proportion $\theta$ of the Bernoulli or binomial distributions. For instance,
when $X_1, \cdots, X_{n}$ are independent and identically distributed (i.i.d.) with Bernoulli($\theta$), i.e., $\Prob_{\theta}(X_i = 1) = 1 - \Prob_{\theta}(X_i = 0) = \theta,$ it is well-known that the maximum likelihood estimator (MLE) of $\theta$ is $\bar X_{n} = (X_1+ \cdots + X_{n})/ n.$ Also  the minimax estimator of $\theta$ under the squared error loss over $\theta \in \Omega = [0,1]$ is
\begin{eqnarray} \label{eqn01}
\delta_{M} = \frac{\sqrt{n}}{\sqrt{n}+1} \bar{X}_{n} + \frac{1}{\sqrt{n} + 1} \frac12,
\end{eqnarray}
%The proof of the minimax property of the estimator in (\ref{eqn01}) is based on the well-known result that a Bayes procedure with constant risk is minimax,
see  \cite{lehmann2006theory}. Here we consider the same minimax estimatior problem of the Bernoulli distribution when $n=1$ but with a twist of the restricted parameter space  $\Omega = \{0 \le \theta \le \eta\}$ for some $0 \le \eta \le 1.$

%When we are teaching the minimax estimator,
The motivating example of our setup can be illustrated through estimating the positive rate $\theta$ of COVID-19 surveillance tests. As an  illustration, suppose that an organization is pretty sure that the COVID-19 positive rate within its employee/staff $\theta \in [0,\eta]$, say, $\eta = 0.2,$
and also runs a surveillance testing program to monitor the current state of the pandemic. Suppose that a random employee/staff is chosen to be sampled at a given day, and how to update the estimated positive rate if the test is positive? What happens if the test is negative? 

%Suppose that the organization tests $n=100$ random employees/staff for COVID-19 surveillance tests during a given period (say daily, weekly or monthly), and observes $5$ positive cases. In this case, the standard sample mean estimator is to estimate the positive rate is $\bar X = \frac{5}{100} = 0.05.$ Meanwhile, the minimiax estimator in (\ref{eqn01}) yields the estimate
%$
%\delta_{M} = \frac{\sqrt{100}}{1 + \sqrt{100}} 0.05 + \frac{1}{1 + \sqrt{100}} \frac12 = 0.0939,
%$
%which is much larger than the sample mean of $5\%.$ However, while the estimator in (\ref{eqn01}) is minimax when the sample space $\Omega = \{0 \le \theta \le 1\},$ it is interesting to ask whether it is still minimax when the sample space $\Omega = \{0 \le \theta \le \eta\}$ for some $0 \le \eta \le 1.$
%
%Meanwhile, the COVID-19 example also brings a new challenge. For instance, suppose that we are pretty sure that the COVID-19 positive rate $\theta \in [0,\eta]$, say, $\eta = 0.2,$ instead of $\theta \in [0,1].$ While the estimator in (\ref{eqn01}) is minimax when the sample space $\Omega = \{0 \le \theta \le 1\},$ it is interesting to ask whether it is still minimax when the sample space $\Omega = \{0 \le \theta \le \eta\}$ for some $0 \le \eta \le 1.$
%This is when the Bayes procedure has non-constant risk function, often when the domain $\Omega$ for the parameter $\theta$ is restricted.

In this note, we provide an elementary mathematical  method to find the minimax estimator for the Bernoulli distribution when $n=1$ and the sample space $\Omega = \{0 \le \theta \le \eta\}$. Our method is to directly solve the optimization method instead  of using the Bayesian device, and thus our proof can be understood by students with diverse backgrounds. 
%We will show that  the estimator in (\ref{eqn01}) is still minimax when $\eta \ge 3/4,$ but we will have a new minimax estimator when $0 \le \eta < 3/4.$ 
%Our proof is based on elementary mathematics tools, as we focus on the $n=1$ case. 
We hope that our note provides more examples of minimax procedures that can be accessed to high school or undergraduate students, thereby enriching the teaching materials on the minimax estimator and enhancing students' understanding and interests.

We should mention that our paper essentially deals with the minimax estimators in restricted parameter spaces, and there are some related existing research for Bernoulli or Binomial distribution in the literature. For instance,  Moors \cite{Moors1985minimax} or  Berry \cite{berry1989bayes}  derived  the minimax estimation when the parameter space $\Omega$ is a symmetric interval, i.e., $\Omega = \{\frac{1}{2}-\eta \le \theta \le \frac{1}{2}+\eta\}$.   Marchand and MacGibbon \cite{marchand2000minimax} considered the parameter space $\Omega = \{ 0 \le \theta \le \eta\}$ as in our note, and derived the minimax estimator for general $n$ case. When $n=1$, their results are equivalent to ours but their methods are too complicated to be accessible to high school, undergraduate or master students. Here, we present a simple method  that can be explained to the students with only elementary mathematics. 

The remainder of this note is as follows. In Section \ref{Sec:02}, we present the minimax estimation problem of Bernoulli proportion in the restricted parameter spaces when the parameter space  $\Omega = \{ 0 \le \theta \le \eta\}$. In Section \ref{Sec:03}, we state our main result on the minimax estimator, and present a rigorous elementary mathematical proof. Section \ref{Sec:04} contains some conclusion marks.

%==========================================================================
\section{Problem Formulation} \label{Sec:02}

%We consider the problem $(P)$ in the following form:\\
% \indent

Suppose that $X$ is a Bernoulli random variable with $P_\theta(X=1) = 1- P_\theta(X=0) = \theta$, and we want to estimate $\theta$ on the basis of $X$ under the squared error loss function $L(\theta,d) =(\theta -d)^2$. In this case, we have $n=1$ observation, and the MLE estimator is to estimate $\theta$ as $\hat\theta_{MLE} = 0$ if $X=0$ and  $=1$ if $X=1.$ Meanwhile, the minimax estimator $\delta_{M}$ in (\ref{eqn01}) is to estimate $\theta$ as $1/4$ if $X=0$ and $3/4$ if $X=1.$

Here we assume further that the parameter $\theta$ belongs to a restricted space $\Omega =\{ 0 \le \theta \le \eta\}$ for some pre-specified constant $\eta,$ and we want to find the minimiax estimator of $\theta$ under the squared loss function. That is, we want to find a procedure $\delta(X)$ that minimizes the maximum risk
\begin{equation} \label{riskf}
\sup_{0 \le \theta \le \eta} \E_{\theta}( \theta - \delta(X))^2
\end{equation}

Since there is only $n=1$ observations, the estimator $\delta(X)$ is completely determined by $\delta(X=0) = a$ and $\delta(X=1) = b.$ Indeed, using the notation of $a$ and $b,$ we have
\begin{eqnarray} \label{eqn03}
\E_{\theta}( \theta - \delta(X))^2 &=& (\theta - a)^{2} \Prob_{\theta}(X= 0 ) + (\theta - b)^{2} \Prob_{\theta}(X= 1 )  \cr
&=& (\theta-a)^2 (1-\theta) + (\theta - b)^2 \theta \cr
&=& (2a - 2b +1)\theta^2 + (b^2 - a^2 -2a )\theta + a^2.
\end{eqnarray}

Thus, the minimax estimator problem in (\ref{riskf}) can be written in the following elementary mathematical form:
\begin{align}
\mbox{ Find two real numbers, $a$ and $b$, that minimize} \cr
\sup_{0 \le \theta \le \eta} \left[ (2a - 2b +1)\theta^2 + (b^2 - a^2 -2a )\theta + a^2 \right]. \label{minmax}
\end{align}
for some $0 \le \eta \le 1.$

%==========================================================================
\section{Our Main Result} \label{Sec:03}

The following theorem presents the solution to the optimization problem in \eqref{minmax}:

\begin{theorem}\label{Theorem2}
When $\Omega = \{ 0 \le \theta \le \eta\},$ the minimax estimator of $\theta$ under (\ref{riskf}), or the optimal solution in (\ref{minmax}), is given by
\[
\delta(X=0) = a^{*} =
	\begin{cases}
        \sqrt{1 - \eta} -(1-\eta) &\quad\text{if $0 \leq \eta \leq \frac{3}{4}$}\\
       \frac{1}{4} &\quad\text{if $\frac{3}{4}< \eta \leq 1$} \\
     \end{cases}
\]
and
\[
\delta(X=1) = b^{*} =  \min \{ \eta, \frac{3}{4} \} = \left\{
                                                        \begin{array}{ll}
                                                          \eta, & \hbox{if $\eta \le \frac34$;} \\
                                                          \frac34, & \hbox{if $\eta > \frac34$.}
                                                        \end{array}
                                                      \right.
\]
%	The corresponding minimax value of  problem $(P)$ is $\min\{\eta^2 -3\eta +2 -2(1-\eta)\sqrt{1-\eta}, \frac{1}{16}\}$.
%	Similarly, we can perform the minimax estimation on the problem $(P)$, where $\Omega = \{\theta \mid \theta \in [\eta,1] , \ 0\leq \eta < 1 \}$.
\end{theorem}

In particular, Theorem \ref{Theorem2} indicates that the estimator $\delta_{M}$ in (\ref{eqn01}) is still minimax for the restricted parameter space $\Omega = \{ 0 \le \theta \le \eta\}$ when $\frac34 \le \eta \le 1,$ but we will have a new form of minimax estimators when $0 \le \eta  \le \frac34.$ For instance,  assume that we are pretty sure that the COVID-19 positive rate of an organization is $\theta \in [0,\eta]$, say, $\eta = 0.2,$ and assume that we randomly test a subject. Then the minimax estimate of the positive rate $\theta$ is $ \sqrt{1 - 0.2} -(1-0.2) \approx 9.44\%$  if the subject's test is negative, and is $20\%$ if the subject's test is positive.

\bigskip
Let us now provide a rigorous proof of Theorem \ref{Theorem2}, which was included in \cite{marchand2000minimax}, but here we present a different proof that directly solves the optimization problem in \eqref{minmax}.

\noindent {\bf Proof of Theorem \ref{Theorem2}}:  It suffices for us to consider $0 \le a \le \eta$ and $0 \le b \le \eta$ when solving the optimization problem in (\ref{minmax}) for $0 \le \eta \le 1.$ This is because it is evident from (\ref{eqn03}) that for all $0 \le \theta \le \eta,$ we have $(\theta - a)^{2} \ge (\theta - \eta)^2$ if $a > \eta$ and $(\theta - a)^{2} \ge (\theta - 0)^2$ if $a < 0.$ Thus, the optimal solution in (\ref{minmax}) must be in the interval $[0,\eta],$ since otherwise we will improve the objective function if we replace them by the endpoints $0$ or $\eta.$

Moreover, note that the function
\begin{equation} \label{eqn05}
f(a,b,\theta) = \left[ (2a - 2b +1)\theta^2 + (b^2 - a^2 -2a )\theta + a^2 \right]
\end{equation}
is a quadratic function with respect to $\theta$, the investigation of its maximum values depends on the sign of the leading coefficient,
and thus we need to split the region of $(a,b)$ into two sub-regions:
\begin{description}
  \item[(A)] When $2a- 2b+1 \ge 0$, the maximum of $f(a,b,\theta)$ is attained at one of two endpoints, $\theta = 0$ or $\theta = \eta.$
  \item[(B)] When $2a- 2b+1 < 0$, $f(a,b,\theta)$ is a concave quadratic function of $\theta$
\end{description}
Our main idea is to find the maximum value that can be minimized in each sub-region, which will yield the global minimax solution.

Let us begin with case (A), or when $a \ge b - \frac{1}{2},$ it suffices to compare two endpoints. In this case, $f(a,b,0) = a^2$ and $f(a,b,\eta) = (1-2b+2a)\eta^2 + (b^2-2a-a^2)\eta +a^2.$ It is evident that $f(a,b,0) \ge f(a,b,\eta)$ if and only if $(1-2b+2a)\eta + (b^2-2a-a^2) \leq 0$, or equivalently, $a^2 + 2(1-\eta) a \ge b^2 - 2 \eta b + \eta.$ This can be further simplified as
\begin{equation} \label{eqn06}
(a + (1-\eta))^2 - (b -\eta)^2 \ge 1-\eta,
%\frac{(a -(1-\eta))^2}{(\sqrt{1-\eta})^2}-\frac{(b-\eta)^2}{(\sqrt{1-\eta})^2} =1,
\end{equation}
or $a \ge \sqrt{ (\eta - b)^2 + (1-\eta)} - (1-\eta),$ since we are only interested in the case when $a > 0.$ A key observation is that  the boundary of (\ref{eqn06}) defines a hyperbolic, which has a vertex $(a,b) = (\sqrt{1-\eta}-(1-\eta), \eta)$ that has the smallest positive $a$ value. Moreover, in case (A), we have  $a \ge b - \frac{1}{2},$ and its boundary line $a = b - \frac{1}{2}$ interests the  hyperbolic boundary of (\ref{eqn06}) at the unique point $(a,b) = (\frac14, \frac34).$ This leads to two subcases, depending on whether $\eta \le \frac34$ or not. This determines whether the unique intersection point is below or above  the vertex of the  hyperbolic curve, or equivalently, whether the line $a = b - \frac{1}{2}$  intersects on the upper or lower part of the hyperbolic curve.

Meanwhile, in case (A), we will need to solve two sub-problems:
\begin{align}
{\bf Problem (A1):}	\min_{0 \le a,b \le \eta} & \ a^2  \label{LskewedP}\\ \nonumber
		s.t. & \ a \ge \sqrt{ (\eta - b)^2 + (1-\eta)} - (1-\eta) \\ \nonumber
		    & \ a \geq b -\frac{1}{2}.
	\end{align}

\begin{align}
{\bf Problem (A2):}	\min_{0 \le a,b \le \eta} & \Big( (1-2b+2a)\eta^2 + (b^2-2a-a^2)\eta +a^2 \Big) \label{RskewedP}\\ \nonumber
		s.t. & \ a \leq \sqrt{ (\eta - b)^2 + (1-\eta)} - (1-\eta) \\ \nonumber
		    & \ a \geq b -\frac{1}{2}.
	\end{align}
In other words, in case (A), we need to investigate two subproblems, (A1) and (A2), under two subcases: one for $0 \le \eta \le \frac34$ and the other for $\frac34 \le \eta \le 1.$

First, we claim that when $0 \le \eta \le \frac34$,  both problems (A1) and (A2) have the same optimal solution:
\begin{eqnarray} \label{eqn09}
a^{*} = \sqrt{1-\eta}-(1-\eta) \qquad \mbox{ and } \qquad b^{*} = \eta.
\end{eqnarray}
To see this, since $0 \le \eta \le \frac34$ and the line $a = b - \frac{1}{2}$  intersects on the upper part of the hyperbolic curve, the objective function $a^2$ in problem (A1) (\ref{LskewedP}) attains its minimum value of $(\sqrt{1-\eta}-(1-\eta))^2$  at the vertex. The proof for problem (A2) in (\ref{RskewedP}) is a little complicated but can follow similar ideas.

To be more specific, in problem (A2), set the objective value  $(1-2b+2a)\eta^2 + (b^2-2a-a^2)\eta +a^2 = \gamma$. This is equivalent to the ellipse:
\begin{eqnarray} \label{eqn10}
	 \frac{\ (a-\eta)^2}{\left(\sqrt{\frac{\gamma}{1-\eta}}\right)^2} + \frac{\left(b-\eta \right)^2}{\left(\sqrt{\frac{\gamma}{\eta}}\right)^2} = 1,
\end{eqnarray}
with a center $(\eta,\eta)$. Since the semi-major axis and semi-minor axis of the ellipse are proportional to $\gamma$, minimizing $\gamma$ is equivalent
to finding the smallest ellipse that intercepts with a point $(a, b)$ in the feasible region of $(a,b)$ in (A2):
\[ %begin{equation} \label{eqn11}
\{(a,b): 0 \le a \le \eta, 0 \le b \le \eta, a \le \sqrt{ (\eta - b)^2 + (1-\eta)} - (1-\eta), a \ge b - \frac12\}.
\] %end{equation}
The smallest ellipse is obtained when it intercepts with the curve $a = \sqrt{ (\eta - b)^2 + (1-\eta)} - (1-\eta)$ at one point. For any $0 < \eta \le 1$, by letting $\gamma = (\sqrt{1-\eta}-(1-\eta))^2$, we have an ellipse that intersects with the curve $a = \sqrt{ (\eta - b)^2 + (1-\eta)} - (1-\eta)$ at the vertex $(\sqrt{1-\eta}-(1-\eta), \eta)$. Hence, when $0\le \eta \le \frac{3}{4}$, the optimal solution is $(\sqrt{1-\eta}-(1-\eta), \eta)$.

%When $\frac{3}{4} < \eta \leq 1$,  $(\sqrt{1-\eta}-(1-\eta), \eta)$ is not in the feasible region. Thus, the optimal solution is the interception $(\frac{1}{4}, \frac{3}{4})$. We illustrate our discussion in Figure 2.

Second, we claim that when $\frac34 \le \eta \le 1$,  both problems (A1) and (A2) have the same optimal solution:
\begin{eqnarray}  \label{eqn11}
a^{*} = \frac14 \qquad \mbox{ and } \qquad b^{*} = \frac34.
\end{eqnarray}
To prove this, note that when  $\frac34 \le \eta \le 1$,  the line $a = b - \frac{1}{2}$  intersects on the lower part of the hyperbolic curve, i.e., the vertex in (\ref{eqn09}) does not satisfy the constraint $a \geq b -\frac{1}{2}$ for problem (A). In particular, for problem (A1), the point with the smallest $a$ value is $(a,b) = (\frac14, \frac34)$, not the vertex, and thus (\ref{eqn11}) is the optimal solution to problem (A1).  When $\frac{3}{4} \leq \eta \leq 1$, the argument to problem (A2) is similar to those in (A1), as the vertex does not belong to the feasible region. In this case, from the shape of the hyperbola, we see that the optimal solution is also given by (\ref{eqn11}),  see Figure 1 and 2 for the illustration of the optimal solutions, depending on the value of $\eta$.

%%Figure 1
\begin{figure}[h]
	
	\centering
	\includegraphics[width=7cm, height=6cm]{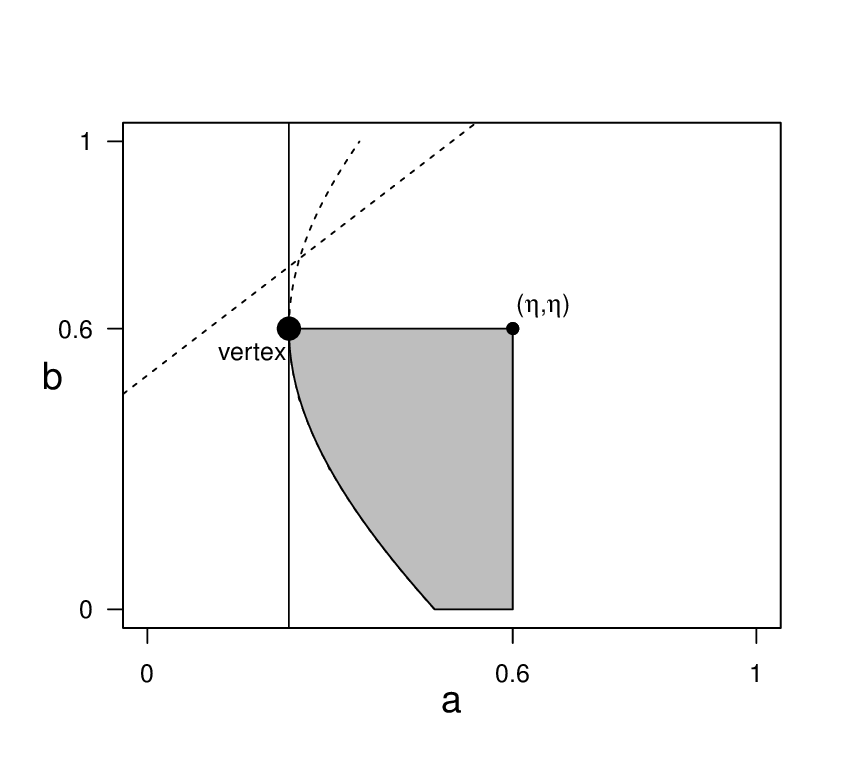}
	\hspace{3 mm}
	\includegraphics[width=7cm, height=6cm]{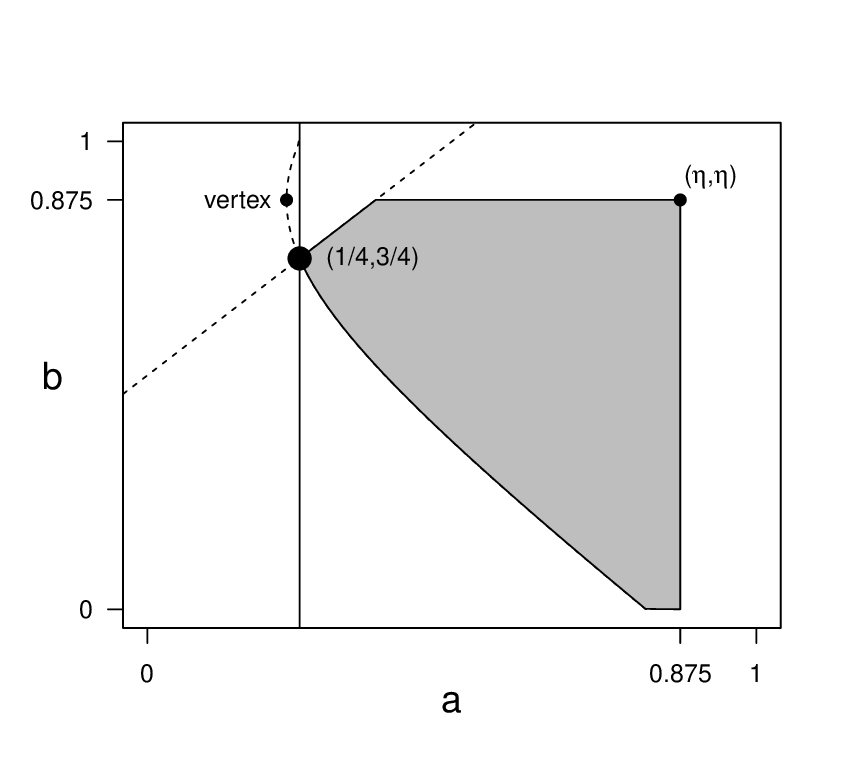}
	
	\caption{The gray area represents the feasible region of problem (A1). The left figure shows that the optimal solution is obtained at the vertex $(\sqrt{1-\eta}-(1-\eta), \eta)$ when $0<\eta \le \frac{3}{4}$. The right figure demonstrates that the optimal solution is $(\frac{1}{4}, \frac{3}{4})$ when $\frac{3}{4} \leq \eta \leq 1$. }
\end{figure}

 %%Figure 2
 \begin{figure}[h]
	
	\centering
	\includegraphics[width=7cm, height=6cm]{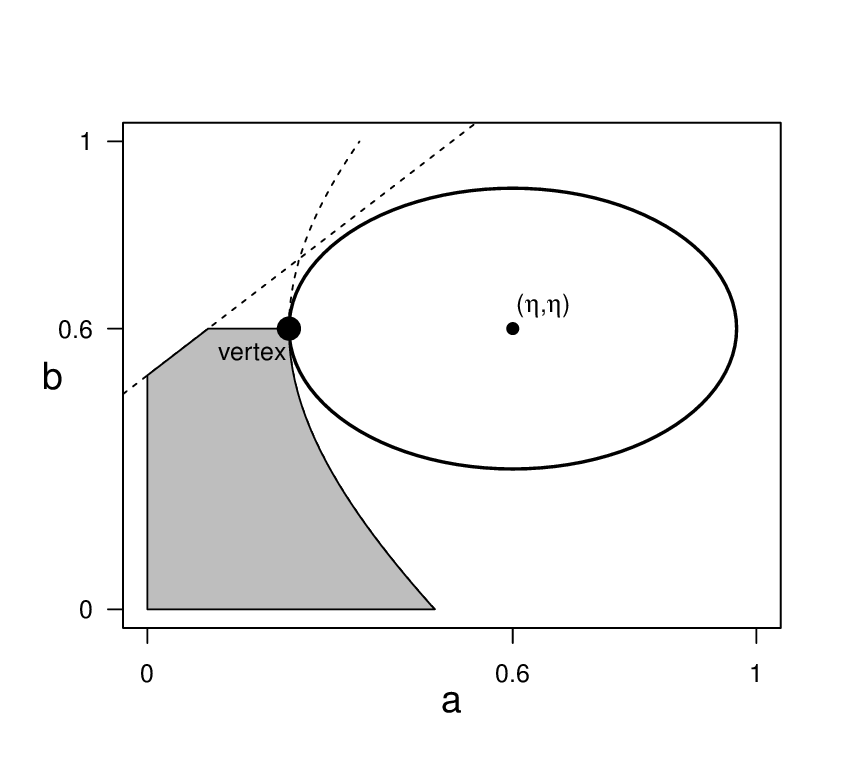}
	\hspace{3 mm}
	\includegraphics[width=7cm, height=6cm]{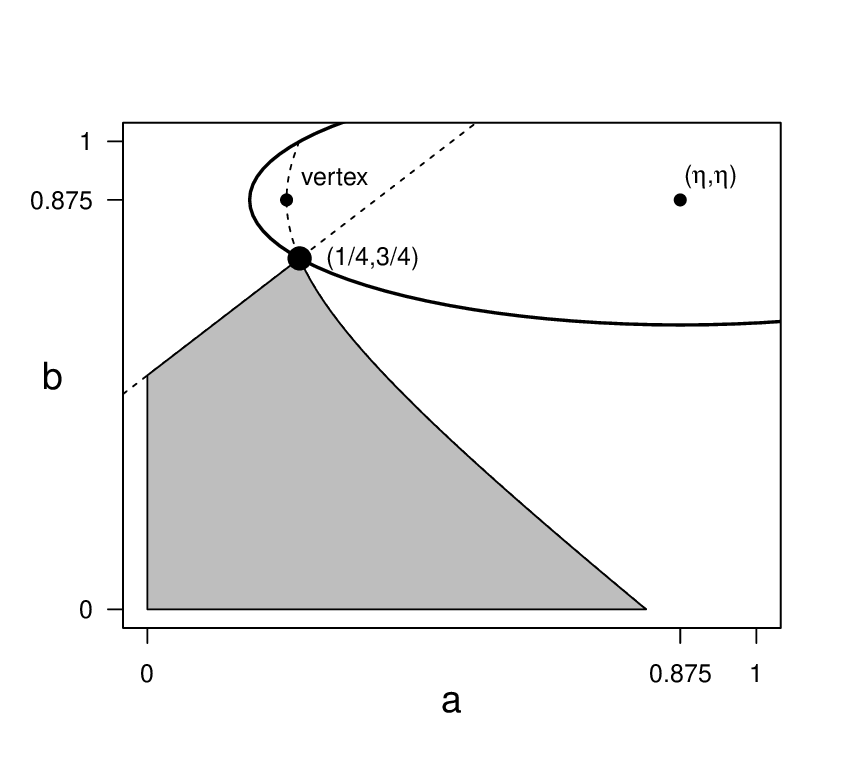}
	
	\caption{The ellipse is centered at $(\eta,\eta)$ and the gray area represents the feasible region of problem (A2). The left figure is
	for the case when $0<\eta \le \frac{3}{4}$ and the right figure is for the case when
	$\frac{3}{4} \leq \eta \leq 1$.}
	% The large black point is the optimal solution, an interception between the smallest ellipse and the curve $a =\sqrt{(b-\eta)^2 +(1-\eta)}-(1-\eta)$.  }
\end{figure}

%% (B) CASE

Next, let us investigate case (B) when $2a- 2b+1 < 0$, i.e., when $a+\frac{1}{2} <b.$ Our main conclusion is that the global minimax solution cannot be obtained in this case. Recall that in the feasible region, we have $0 \leq a \leq \eta$ and $0 \leq b \leq \eta.$ Thus the relation $a+\frac{1}{2} <b$ cannot hold if $\eta \le \frac12.$  Thus, it suffices to investigate case (B) when $\frac{1}{2} < \eta \leq 1$.
We claim that the maximum value in case (B) will be larger those in case (A) when $\frac{1}{2} < \eta \leq 1$. To prove this, consider a specific value $\theta = \frac{1}{2} \in [0 , \eta]$. Then we have the following inequalities:
	  \begin{align*}
 	 f(a,b,\frac{1}{2}) &= \frac{1}{4}(1-2b+2a) + \frac{1}{2}(b^2-2a-a^2)+a^2\\
%	 &= \frac{1}{2}\left(a^2+b^2-a-b+\frac{1}{2}\right)\\
	 &= \frac{1}{2}\left\{\left(a-\frac{1}{2}\right)^2+\left(b-\frac{1}{2}\right)^2\right\} \\
	 &> \frac{1}{2}\left\{ \left(a-\frac{1}{2}\right)^2+a^2\right\} \\
%	 & = \frac{1}{2}\left(2a^2-a+\frac{1}{4}\right) \\
%	 & = a^2-\frac{1}{2}a + \frac{1}{8} \\
	 &= \left(a- \frac{1}{4}\right)^2 + \frac{1}{16} \\
	 &\geq \frac{1}{16}.
   	\end{align*}
where the first inequality follows from the assumption of $a+\frac{1}{2} <b$ for Case (B). On the other hand, recall that the maximum values of case (A) attains the minimum value of $(\sqrt{1-\eta}-(1-\eta))^2$ and $\frac{1}{16}$, respectively, depending on whether $\eta < \frac34$ or not. Both of these maximum values are less than or equal to $\frac{1}{16}$. This implies that a minimax estimator for case (B) cannot be a minimax estimator for problem \eqref{minmax}. As a result, we exclude case (B) and conclude that the minimax estimator in case (A) is the minimax estimator of the problem \eqref{minmax}. This proves the theorem.
\hfill \QED

%=======================================================	
	 \section{Conclusion} \label{Sec:04}

In the statistical inference or statistical decision theory course, the restricted minimax problem \eqref{minmax} would serve as a useful advanced exercise for illustrating minimax estimation. Indeed, this problem have been discussed in our advanced undergraduate level or first-year-graduate-level statistics course to students with diverse backgrounds in mathematics, engineering, and computer sciences. To find the minimax estimators for the problem \eqref{minmax}, we have used the geometrical interpretation of hyperbola, ellipse and  convex  functions. Since our method is simple and different from the standard Bayesian approach, it can be illustrated to students with diverse backgrounds. Moreover, our methodology motivates students to tackle statistical problems with different perspectives, and stimulates their interests on statistics.

\end{document}